\documentclass{article}
\usepackage{graphicx}

\usepackage{amsfonts}
\usepackage{amsmath}
\usepackage{amsthm}

\title{Convergence in inhomogeneous consensus processes with positive diagonals}

\author{
Jan Lorenz, \\ Universit{\"a}t Bremen, Fachbereich Mathematik und
Informatik\\ Bibliothekstra{\ss}e, 28359 Bremen, Germany \\
\texttt{math@janlo.de}}
\date{August 31, 2006}

\newcommand{\n}{\underline n}
\newcommand{\N}{\mathbb N}
\newcommand{\R}{\mathbb R}
\newcommand{\Ci}{\mathcal{I}}
\newcommand{\Cj}{\mathcal{J}}
\newcommand{\spann}{\mathrm{span}}
\newcommand{\eig}{\mathrm{eig}}
\newcommand{\one}{\mathbf{1}}

\newtheorem{theorem}{Theorem}

\newtheorem{proposition}[theorem]{Proposition}

\begin{document}

\maketitle

\begin{abstract}
We present a results about convergence of products of
row-stochastic matrices which are infinite to the left and all
have positive diagonals. This is regarded as in inhomogeneous
consensus process where confidence weights may change in every
time step but where each agent has a little bit of self
confidence. The positive diagonal leads to a fixed zero pattern in
certain subproducts of the infinite product.

We discuss the use of the joint spectral radius on the set of the
evolving subproducts and conditions on the subprodutcs to ensure
convergence of parts of the infinite product to fixed
rank-1-matrices on the diagonal.

If the positive minimum of each matrix is uniformly bounded from
below the boundedness of the length of intercommunication
intervals is important to ensure convergence. We present a small
improvement. A slow increase as quick as $\log(\log(t))$ in the
length of intercommunication intervals is acceptable.
\end{abstract}


\section{Introduction}

Consider $n$ persons that discuss an issue which can be
represented as a real number. Assume further that the persons
revise their opinions if they hear the opinions of others. Each
person finds his new opinion as a weighted arithmetic mean of the
opinions of others. This model of opinion dynamics has been
analyzed for the possibilities of consensus by DeGroot
\cite{DeGroot1974}. If these weights change over time we have an
inhomogeneous consensus process.

While the homogeneous process has strong similarities with a
homogeneous Markov chain, things get different when inhomogeneity
comes in. While a consensus process relies on row-stochastic
matrices multiplied from the left, a Markov process relies on
row-stochastic matrices multiplied from the right. And infinity to
the right is not the same as infinity to the left.

Consensus processes are only briefly touched in the context of
Markov chains \cite{Hartfiel1998}. Besides the early approaches of
opinion dynamics \cite{DeGroot1974, Lehrer1981} some results have
been made in the context of decentralized computation
\cite{Tsitsiklis1984}. Consensus processes fit in the framework of
questions about sets of matrices which have the left convergence
property 'LCP' \cite{Daubechies1992,Hartfiel2002}, which is 'RCP'
for transposed matrices.

Recently, there have been independent works that study consensus
processes and the underlying matrix-products in the context of
opinion dynamics \cite{Krause2000, Hegselmann2002, Lorenz2003b},
multi agent systems where agents try to coordinate
\cite{Jadbabaie2003, Moreau2005} and flocking where birds or
robots try to find agreement about their headings
\cite{Blondel2005a,Hendrickx2005}.

In \cite{Blondel2005a,Theys2005} there have been the first
attempts to make the concept of the joint spectral radius work on
consensus processes.

In this paper we want to analyze the structure that positive
diagonals deliver in inhomogeneous consensus processes and extend
the basic idea of \cite{Blondel2005a,Theys2005}. But a result on
convergence is only reachable with further assumptions on
matrices. In the end we will derive a small improvement on
acceptable growth of the length of intercommunication intervals.

\section{Consensus Processes}

For $n\in\N$ we define $\n := \{1,\dots,n\}$.

Let $A(0), A(1), \dots$ be a sequence of square row-stochastic
matrices of size $n\times n$.

For natural numbers $s<t$ we define a \emph{forward accumulation}
$A(s,t) := A(s)\dots A(t-1)$ and a \emph{backward accumulation}
$A(t,s) := A(t-1)\dots A(s)$. Thus $A(s,s+1)=A(s+1,s)=A(s)$ and
$A(s,s)$ is the identity.

Let $x(0)$ be a real column vector of opinions and $x_i(0)$ stands
for the initial opinion of person $i$. The sequence of vectors
$x(t) = A(t,0)x(0)$ is an \emph{inhomogeneous consensus process}
and $a(t)_{ij}$ stands for a confidence weight person $i$ gives to
the opinion of agent $j$ at time step $t$. In this context $A(t)$
is called a \emph{confidence matrix}.

To understand the convergence behavior of inhomogeneous consensus
processes the infinite product $A(\infty,0)$ is of interest.

In this paper we focus on confidence matrices with positive
diagonals. Thus, we regard processes where persons always have a
little bit of self-confidence.

A row-stochastic matrix $K$ which has rank 1 and thus equal rows
is called a \emph{consensus matrix} because for a real vector $x$
it holds that $Kx$ is a vector with equal entries and thus
represents consensus among persons in a consensus process. Suppose
that $A(t):=K$ is a consensus matrix. It is easy to see that for
all $u\geq t$ it holds for the backward accumulation that $A(u,0)
= K$.  (For the infinite forward accumulation $A(0,\infty)$ it
only holds that $A(0,u)$ is a consensus matrix but may change with
$u$.) In the following we will point out that there is also a
tendency of convergence to consensus matrices.

In the next section we will see that the positive diagonal
together with the Gantmacher's canonical form of nonnegative
matrices \cite{Gantmacher1959_2} will give us a good overview on
the zero and positivity structure of the processes.

In section \ref{sec:conv} we go on with a convergence theorem that
is built on this structure and conclude in section \ref{sec:disc}
with a small improvement and discussion on how to fulfill the
conditions of the theorem.

\section{The positive diagonal} \label{sec:posdiag}

We regard two nonnegative matrices $A,B$ to be of the same
\emph{type} $A \sim B$ if $a_{ij} > 0 \Leftrightarrow b_{ij}>0$.
Thus, if their zero-patterns are equal. All matrices of the same
type have the same Gantmacher form, which block structure we will
outline now.

Let $A$ be a nonnegative matrix with a positive diagonal. For
indices $i,j\in\n$ we say that there is a \emph{path} $i
\rightarrow j$ if there is a sequence of indices
$i=i_1,\dots,i_k=j$ such that for all $l\in\underline{k-1}$ it
holds $a_{i_l,i_{l+1}}>0$. We say $i,j \in\n$ \emph{communicate}
if $i \rightarrow j$ and $j \rightarrow i$, thus $i
\leftrightarrow j$. In our case with positive diagonals there is
always a path from an index to itself, which we call
\emph{self-communicating} and thus $'\leftrightarrow'$ is an
equivalence relation. An index $i\in\n$ is called \emph{essential}
if for every $j\in\n$ with $i \rightarrow j$ it holds $j
\rightarrow i$. An index is called inessential if it is not
essential.

Obviously, $\n$ divides into disjoint self-communicating
equivalence classes of indices $\Ci_1,\dots,\Ci_p$. Thus, in one
class all indices communicate and do not communicate with other
indices. The terms essential and inessential thus extend naturally
to classes. We define $n_1 := \#\Ci_p,\dots,n_p := \#\Ci_p$.

If we renumber indices with first counting the essential classes
and second the inessential classes with a class $\Ci$ before a
class $\Cj$ if $\Cj\rightarrow\Ci$ then we can bring every
row-stochastic matrix $A$ to the \emph{Gantmacher form}
\cite{Gantmacher1959_2}
\begin{equation} \label{gf}
\begin{bmatrix}
  A_1 &  &  & &  & 0 \\
   & \ddots &  &  &  &  \\
  0 &  & A_g &  &  &  \\
  A_{g+1,1} & \dots & A_{g+1,g} & A_{g+1} &  &  \\
  \vdots &  & \vdots & \vdots & \ddots &  \\
  A_{p,1} & \dots & A_{p,g} & A_{p,g+1} & \dots & A_p \\
\end{bmatrix}
\end{equation}
by simultaneous row and column permutations. The \emph{diagonal
Gantmacher blocks} $A_1,\dots,A_p$ in (\ref{gf}) are square ($n_1
\times n_1,\dots,n_p \times n_p$) and irreducible. Irreducibility
induces primitivity in the case of a positive diagonal. For the
\emph{nondiagonal Gantmacher blocks} $A_{k,l}$ with
$k=g+1,\dots,p$ and $l=1,\dots,k-1$ it holds that for every
$k\in\{g+1,\dots,p\}$ at least one block of
$A_{k,1},\dots,A_{k,k-1}$ contains at least one positive entry.

The spectrum of $A$ is the union of the spectra of all the
diagonal Gantmacher blocks.

The following proposition shows that an infinite backward or
forward accumulation of nonnegative matrices can be divided after
a certain time step into subaccumulations with a common Gantmacher
form.

\begin{proposition} \label{prop1}
Let $(A(t))_{t\in\N}$ be a sequence of nonnegative matrices with
positive diagonals. Then for the backward accumulation there
exists a sequence of natural numbers $0 < t_0 < t_1 < \dots$ such
that for all $i\in\N$ it holds
\begin{equation}
A(t_{i+1},t_{i}) \sim A(t_1,t_0).
\end{equation}
Thus, $A(t_{i+1},t_{i})$ can be brought to the same Gantmacher
form for all $i\in\N$. Further on, all Gantmacher diagonal blocks
are positive and all nondiagonal Gantmacher-Blocks are either
positive or zero.
\end{proposition}

\begin{proof}
(In sketch, for more details see \cite{Lorenz2005}.)

The proof works with a double monotonic argument on the positivity
of entries: While more and more (or exactly the same) positive
entries appear in $A(t,0)$ monotonously increasing with rising
$t$, we reach a maximum at $t^\ast_0$. We cut $A(t^\ast_0,0)$ of
and find $t^\ast_1$ when $A(t,t^\ast_0)$ reaches maximal
positivity again with rising $t$. We go on like this and get the
sequence $(A(t^\ast_{i+1},t^\ast_i))_{i\in\N}$. Obviously, less
and less (or exactly the same) positive entries appear
monotonously decreasing with rising $i$ and we reach a minimum at
$k$. We relabel $t_j := t^\ast_{k+j}$ and thus have the desired
sequence $(t_i)_{i\in\N}$ with $A(t_{i+1},t_i)$ having the same
zero-pattern.

Positivity of Gantmacher blocks follows for all blocks
$A(t_{i+1},t_i)_{[\Cj,\Ci]}$ where we have a path $\Cj \rightarrow
\Ci$. If we have such a path, then there is a path from each index
in $\Cj$ to each index in $\Ci$ and thus every entry must be
positive in a long enough accumulation. Thus, the block has to be
positive already, otherwise $(t_i)_{i\in}$ is chosen wrong.

To prove the result for forward accumulations, we can use the same
arguments.
\end{proof}

Let us consider now a sequence of row-stochastic matrices
$(A(t))_{t\in\N}$ and their infinite backward products $A(t,0)$
with $t\to\infty$. Thus, we face a consensus process where agents
may change their confidence weights in every time step.

Form proposition \ref{prop1} we get the existence of a sequence of
time step $(t_i)_{i\in\N}$ such that all $A(t_{i+1},t_i)$ have the
same Gantmacher form with positive Gantmacher diagonal blocks. So,
the Gantmacher structure represents, that agents find a stable
confidence structure. There evolve $g \geq 1$ groups where every
agents trust everyone else internally (but maybe indirectly) and
no one outside; this repeats for all the time. And there evolve
inessential confidence groups in which agents trust each other
internal but which also have trust chains to one or more of the
$g$ essential groups.

Unfortunately, nothing can be said about the distances
$t_{i+1}-t_i$.

\section{The joint spectral radius}

We regard a sequence of row-stochastic matrices with positive
diagonals $(A(t))_{t\in\N}$, take the sequence of time steps of
proposition \ref{prop1} and abbreviate $A(i) := A(t_{i+1},t_i)$.
Further on, the $A_k(i),A_{k,j}(i)$ are the respective Gantmacher
blocks of $A(i)$. So, $\Sigma :=\{A(i) \,|\, i\in\N\}$ is a set of
matrices with the same Gantmacher form, which joint spectral
radius can be studied.

The spectral radius of a matrix $A$ is $\rho(A) := \{|\lambda|
\,\, \lambda \hbox{ is eigenvalue of } A\}$ and represents the
growth rate of the matrix norm of $A^i$. The joint spectral radius
\cite{Daubechies1992} of a set of square matrices $\mathcal{M}$ is
\[
\hat\rho(\mathcal{M}) := \limsup_{k\to\infty}
\sup_{A(i_1),\dots,A(i_k) \in \mathcal{M}} ||A(i_1)\dots A(i_k)
||^\frac{1}{k}
\]
and represents the maximal growth rate of arbitrary products of
matrices from $\mathcal{M}$.

In our setting for all $i \in \N$ it holds $\rho(A(i)) = 1$ and
due to the fact that every product of $\Sigma$ is row-stochastic
it holds $\hat\rho(\Sigma) = 1$, too. But we can do a joint
transformation of all matrices in $\Sigma$ which leads us to a
situation where the joint spectral radius is more interesting.

Let us consider the $k$-th Gantmacher diagonal block $A_k$ for the
essential class $\Ci_k$ in an arbitrary accumulation
$A_k(t_{i+1},t_i) =: A_k(i)$ ($k\in\underline{g}$, $i\in\N$).
$A(i)_k$ is positive and row-stochastic. Thus, it has the unique
maximal eigenvalue $1$ for the eigenvector $\one$. ($\one$ is the
vector with only one-entries of the appropriate length given
through the context). And there are no other eigenvalues with
absolute value one.

According to an idea outlined in \cite{Blondel2005a,Theys2005} we
can make a transformation $P_kA(i)_kP_k^T =: A'(i)$ such that the
spectrum stays the same but with eigenvalue 1 removed. For this
$P_k$ is a $(n_k-1) \times n$ matrix which rows build an
orthogonal basis of the orthogonal complement to $\spann\{\one\}$.
(This can be normalized vectors with two nonzero entries which
have the same absolute value and different signs.)

Thus, $A'(i)$ is $(n_k-1)\times (n_k-1)$. To see that the spectrum
of $A'(i)_k$ is the spectrum of $A(i)_k$ without 1 consider an
eigenvalue $\lambda \neq 1$ and one of its eigenvectors $x$. Then
$y := P_kx$ is not zero and an eigenvector of $A'(i)_k$ for the
eigenvalue $\lambda$. ($A(i)_kx=\lambda x \Rightarrow
P_kA(i)_kP_k^Ty = \lambda y \Rightarrow A'(i)_ky=\lambda y$.)

Obviously, all the matrices $A(i)$ have 1 as eigenvalue $g$ times
with a $g$-dimensional eigenspace
\begin{equation}\label{eq:span}
\eig(A(i),1) = \spann\{
\left[%
\begin{array}{c}
  \one \\
  0 \\
  \vdots \\
  0 \\
  \ast \\
\end{array}%
\right],
\left[
\begin{array}{c}
  0 \\
  \one \\
  \vdots \\
  0 \\
  \ast \\
\end{array}
\right], \dots,
\left[
\begin{array}{c}
  0 \\
  \vdots \\
  0 \\
  \one \\
  \ast \\
\end{array}
\right] \}
\end{equation}

The $\ast$-parts are not necessary equal for all $i\in\N$, but it
is clear that $\one$ is in, thus the $\ast$-parts sum up to
multiple of $\one$.

Nevertheless, we can generalize the transformation idea of
\cite{Blondel2005a,Theys2005} to our setting. We define the
$(n-g)\times n$ matrix
\[P := \left[%
\begin{array}{cccc}
  P_1 &  &  & 0 \\
   & \ddots &  &  \\
   &  & P_g &  \\
  0 &  &  & E \\
\end{array}%
\right]\] where $E$ is the unit matrix of size
$n_{g+1}+\dots+n_p$. Notice that the blocks are not square and
thus not diagonal. Now, it holds $P A(t_{i+1},t_i) P^T =$
\[
\left[%
\begin{array}{cccccc}
  A'_1(i) &  &  & & & 0 \\
   & \ddots &  & & &  \\
  0  &  & A'_g(i) &  & & \\
  A_{g+1,1}(i)P_1^T & \dots & A_{g+1,g}(i)P_g^T & A_{g+1}(i) & & \\
  \vdots &  & \vdots &  & \ddots & \\
  A_{p,1}(i)P_1^T & \dots & A_{p,g}(i)P_g^T & A_{p,g+1}(i) & \dots & A_{p}(i) \\
\end{array}%
\right] =: A'(i)
\]

Now we can study the joint spectral radius of $\Sigma' :=
{PA(i)P^T | i\in\N}$. If we had $\hat\rho(\Sigma') < 1$ this would
imply that $A(t,0)x(0)$ would converge in the entries of indices
$\Ci_1 \cup \dots \cup \Ci_g$ to a vector in
\[
\spann\{\left[%
\begin{array}{c}
  \one \\
  0 \\
  \vdots \\
  0 \\
\end{array}%
\right], \left[
\begin{array}{c}
  0 \\
  \one \\
  \vdots \\
  0 \\
\end{array}
\right], \dots, \left[
\begin{array}{c}
  0 \\
  \vdots \\
  0 \\
  \one \\
\end{array}%
\right]\}.
\]

It holds for the spectral radii that
$\rho(A'_1(i))<1,\dots,\rho(A'_g(i))<1$, due to the the fact that
$A_1(i),\dots,A_g(i)$ where positive and thus had no other
eigenvalues of absolute value one. Further on, the spectral radii
of $A_{g+1}(i),\dots,A_p(i)$ are less than one because for
$l\in\{g+1,\dots,p\}$ it holds $\rho(A_l(i)) \leq ||{A_l(i)}|| <
1$. The second inequality holds due to the fact that all row sums
in $A_l(i)$ are less than one. ($||A|| := \max_i\sum_j |a_{ij}|$
in this case.)

Thus, it holds $\rho(A'(i)) < 1$ for all $i\in\N$. But
unfortunately this does not imply $\hat\rho(\Sigma')<1$
\cite{Blondel2000}. Thus, more assumptions must be made to reach a
partial convergence result. This is subject to the next section,
where we use concepts of ergodicity.

\section{Convergence} \label{sec:conv}

We define the \emph{coefficient of ergodicity} of a row-stochastic
matrix $A$ according to \cite{Hartfiel1998} as
\[\tau(A) := 1 - \min_{i,j \in \n} \sum_{k=1}^n \min\{a_{ik}, a_{jk}\}.\]

The coefficient of ergodicity of a row-stochastic matrix can only
be zero, if all rows are equal, thus if it is a consensus matrix.

The coefficient of ergodicity is submultiplicative (see
\cite{Hartfiel1998}) for row-stochastic matrices $A_0,\dots,A_i$
\begin{equation}\label{eq:submulterg}
\tau(A_i \cdots A_1 A_0) \leq \tau(A_i) \cdots \tau(A_1)
\tau(A_0).
\end{equation}

If $\lim_{t\to\infty} \tau(A(0,t)) = 0$ we say that $A(0,t)$ is
\emph{weakly ergodic}. Weakly ergodic means that the $A(0,t)$ gets
closer and closer to the set of consensus matrices and thus the
Markov process gets totally independent of the initial
distribution.

For $M \subset \R_{\geq 0}$ we define $\min^{+}M$ as the smallest
positive element of $M$. For a stochastic matrix $A$  we define
$\min^{+} A := \min^{+}_{i,j \in\n} a_{ij}$. We call $\min^{+}$
the \emph{positive minimum}.

For the positive minimum of a set of row-stochastic matrices
$A_0,\dots,A_i$ it holds
\begin{equation}\label{eq:posmin}
\min\mbox{}^{+}(A_i\cdots A_0) \geq \min\mbox{}^{+}A_i\cdots
\min\mbox{}^{+}A_0.
\end{equation}

\begin{theorem} \label{theorem}
Let $(A(t))_{t\in\N}$ be a sequence of row-stochastic matrices
with positive diagonals, $0<t_0<t_1<\dots$ be the sequence of time
steps defined by proposition \ref{prop1}, $\Ci_1,\dots,\Ci_g$ be
the essential and $\Cj$ be the union of all inessential classes of
$A(t_1,t_0)$.

If for all $i\in \N$ it holds $\min^{+} A(t_{i+1},t_i) \geq
\delta_i$ and $\sum_{i=1}^\infty \delta_i = \infty$, then
\[
\lim_{t\to\infty} A(t,0) = \left[\begin{array}{ccc|c}
  K_1 &  & 0 & 0 \\
   & \ddots &  & \vdots \\
  0 &  & K_g & 0 \\ \hline
   & \mathrm{not\ converging} &  & 0 \\
\end{array}\right] A(t_0,0)
\]
where $K_1,\dots,K_g$ are consensus matrices. (The matrices have
to be sorted by simultaneous row and column permutations according
to $\Ci_1,\dots,\Ci_g,\Cj$.)
\end{theorem}

\begin{proof}
The interesting blocks are the diagonal blocks. It is easy to see
due to the lower block triangular Gantmacher form of
$A(t_{i+1},t_i)$ for all $i\in\N$, that all diagonal blocks only
interfere with themselves when matrices are multiplied.

Let us regard the essential class $\Ci_k$ and abbreviate $A_i :=
A(t_{i+1},t_i)_{[\Ci_k,\Ci_k]}$.

We show that the minimal entry in a column $j$ of a row-stochastic
matrix $B$ cannot sink when multiplied from the right with another
row-stochastic matrix $A$,
\[\min_{i\in\n} (AB)_{ij} = \min_{i\in\n} \sum_{k=1}^n a_{ik}b_{kj} \geq \min_{i\in\n} b_{ij}.\]
Thus, the minimum of entries in column $j$ of the product
$A_i\cdots A_0$ is monotonously increasing with rising $i\in\N$.
With similar arguments it follows that the maximum of entries in
column $j$ of the product $A_i\cdots A_0$ is monotonously
decreasing with rising $i\in\N$.

Further on, it holds due to (\ref{eq:submulterg}) and the
definition of the coefficient of ergodicity that
\[\lim_{i\to\infty}\tau(A_i\dots A_1 A_0) \leq
\prod_{i=1}^{\infty}\tau(A_i) = \prod_{i=1}^{\infty}(1-\delta_i)
\leq \prod_{i=1}^{\infty}e^{-\delta_i} =
e^{-\sum_{i=1}^{\infty}\delta_i} = 0.
\]
The maximal distance of rows shrinks to zero. Both arguments
together imply that $\lim_{i\to\infty}(A_i\dots A_1 A_0)$ is a
consensus matrix which we call $K_k$.

Now it remains to show that the $[\Cj,\Cj]$-diagonal block of the
inessential classes converges to zero.

Let us define $||\cdot||$ as the row-sum-norm for matrices. It
holds $||A_{[\Cj,\Cj]}(t_{i+1},t_{i})|| \leq (1-\delta_i)$ and
thus like above it holds
\[||A_{[\Cj,\Cj]}(\infty,t_0)|| \leq \prod_{i=1}^\infty ||A_{[\Cj,\Cj]}(t_{i+1},t_i)|| \leq \prod_{i=1}^\infty (1-\delta_i) \leq = 0.\]
This proves that $\lim_{t\to\infty}A_{[\Cj,\Cj]}(t,0) = 0$.
\end{proof}

An inhomogeneous consensus process $A(t,0)x(0)$ with persons who
have some self-confidence stabilizes (under weak conditions) such
that we have $g$ consensual subgroups (the essential classes)
which have internal consensus, while all other persons (the
inessential indices) may hop still around building opinions as
convex combinations of the values reached in the consensual
groups.

\section{Discussion on conditions for $\min^{+}A(t_{i+1},t_i) \geq
\delta_i$}\label{sec:disc}

One thing where theorem \ref{theorem} stays unspecific is that it
demands lower bounds for the positive minimum of the accumulations
$A(t_{i+1},t_i)$. But, what properties of the single matrices may
ensure the assumption $\min^{+}A(t_{i+1},t_i) \geq \delta_i$ with
$\sum \delta_i = \infty$?

The first idea would be to assume a \emph{uniform lower bound for
the positive minimum} $\delta < \min^{+}A(t)$ for all $t$. But
this is not enough.

Recent independent research
\cite{Moreau2005,Hendrickx2005,Lorenz2003b} has shown that either
 \emph{bounded intercommunication intervals} ($t_{i+1}-t_i < N$
for all $i\in\N$) or \emph{type-symmetry} ($A \sim A^T$) of all
matrices $A(t)$ can be assumed additional to the uniform lower
bound for the positive minimum to ensure the assumptions of
theorem \ref{theorem}. But improvements are possible.

\subparagraph{Bounded intercommunication intervals} Let us regard
$\delta < \min^{+}A(t)$ for all $t\in\N$. If $t_{i+1}-t_i \leq N$
it holds by (\ref{eq:posmin}) that $\min^{+}A(t_{i+1},t_i) \geq
\delta^{N}$ and thus $\sum_{i=0}^\infty \delta^N = \infty$ and
thus theorem \ref{theorem} holds. But $t_{i+1}-t_i$ may slightly
rise as the next two propositions show.

\begin{proposition}
Let $0<\delta<1$ and $a \in \R_{>0}$ then
\begin{equation}
\sum_{n=1}^\infty \delta^{a\log(n)} < \infty \Longleftrightarrow
\delta < e^{-1}.
\end{equation}
\end{proposition}

\begin{proof}
We can use the integral test for the series $\sum_{n=1}^\infty
\delta^{a\log(n)}$ because $f(x) := \delta^{a\log(x)}$ is positive
and monotonously decreasing on $[1,\infty[$.

With substitution $y = \log(x)$ (thus $dx=e^{y}dy$) it holds
\begin{eqnarray*}
\int_1^\infty \delta^{a\log(x)}dx & = & \int_1^\infty
e^{a\log(\delta)\log(x)}dx = \int_1^\infty
e^{a\log(\delta)y}e^{y}dy \\
& = & \int_1^\infty e^{ay(\log(\delta)+1)}dy
\end{eqnarray*}
The integral is finite if and only if $\log(\delta)+1 < 0$ and
thus if $\delta < e^{-1}$.
\end{proof}

\begin{proposition}
Let $0<\delta<1$ and $a \in \R_{>0}$ then
\begin{equation}
\sum_{n=3}^\infty \delta^{a\log(\log(n))} = \infty.
\end{equation}
\end{proposition}

\begin{proof}
We can use the integral test for the series $\sum_{n=1}^\infty
\delta^{a\log(\log(n))}$ because $f(x) := \delta^{a\log(\log(x))}$
is positive and monotonously decreasing on $[3,\infty[$.

With substitution $y = \log(\log(x))$ (thus $dx=e^{(y+e^y)}dy$) it
holds
\begin{eqnarray*}
\int_3^\infty \delta^{a\log(\log(x))}dx & = & \int_1^\infty
e^{a\log(\delta)\log(\log(x))}dx = \int_1^\infty
e^{a\log(\delta)y}e^{y+e^y}dy \\
& = & \int_1^\infty e^{ay(\log(\delta)+1) + e^y}dy
\end{eqnarray*}
The integral diverges because  $ay(\log(\delta)+1) + e^y
\longrightarrow \infty$ as $y\to\infty$.
\end{proof}

Thus, assuming $\min^{+}A(t) > \delta > 0$ for all $t\in\N$ we can
allow a slow growing of $t_{i+1} - t_i$ to fulfill the assumptions
of theorem \ref{theorem}. Acceptable is a growing as quick as
$\log(\log(i))$. If $t_{i+1} - t_i$ grows as $\log(i)$ then it
must hold $\delta > e^{-1} > \frac{1}{3}$. This can only hold if
each row of $A(t)$ contains only two positive entries (due to
row-stochasticity).

\section{Conclusion}

We pointed out the convergence of the zero patterns of
accumulations in inhomogeneous consensus processes with positive
diagonals. It leads to a stable Gantmacher form on accumulations.
We then extended an idea of \cite{Blondel2005a,Theys2005} to a
potential use of the joint spectral radius for the convergence of
inhomogeneous consensus processes but saw that further assumptions
are necessary to reach a convergence result. For this we switched
back to the concept of shrinking coefficients of ergodicity and
could reach a small improvement of former results.

Perhaps the combination of both approaches may lead to a full
characterization of inhomogeneous consensus processes with respect
to convergence and conditions for consensus.

\subsection*{Acknowledgement}
I thank Dirk A. Lorenz for calculus hints.

\bibliographystyle{unsrt}
\bibliography{../../lit}

\end{document}